\documentclass[12pt,a4paper]{amsart} 

\usepackage{hyperref} 
\usepackage{amssymb}
\usepackage{amsrefs}
\usepackage{amsthm} 
\usepackage[all]{xy}

\allowdisplaybreaks

\newcommand{\bd}{\mathbf{d}}
\newcommand{\be}{\mathbf{e}}
\newcommand{\bh}{\mathbf{h}}
\newcommand{\bm}{\mathbf{m}}
\newcommand{\bP}{\mathbf{P}}
\newcommand{\bQ}{\mathbf{Q}}
\newcommand{\bR}{\mathbf{R}}
\newcommand{\blambda}{\mbox{\boldmath $\lambda$}}

\newcommand{\fZ}{\mathfrak{Z}}

\newcommand{\bbA}{\mathbb{A}}
\newcommand{\bbM}{\mathbb{M}}
\newcommand{\bbN}{\mathbb{N}}
\newcommand{\bbP}{\mathbb{P}}
\newcommand{\bbZ}{\mathbb{Z}}

\newcommand{\calC}{\mathcal{C}}
\newcommand{\calO}{\mathcal{O}}
\newcommand{\calP}{\mathcal{P}}
\newcommand{\calQ}{\mathcal{Q}}
\newcommand{\calR}{\mathcal{R}}
\newcommand{\calX}{\mathcal{X}}
\newcommand{\calY}{\mathcal{Y}}
\newcommand{\calZ}{\mathcal{Z}}

\let\mod=\undefined
\DeclareMathOperator{\GL}{GL} %
\DeclareMathOperator{\id}{id} %
\DeclareMathOperator{\pd}{pd} %
\DeclareMathOperator{\SI}{SI} %
\DeclareMathOperator{\SL}{SL} %
\DeclareMathOperator{\Ext}{Ext} %
\DeclareMathOperator{\Hom}{Hom} %
\DeclareMathOperator{\mod}{mod} %
\DeclareMathOperator{\Reg}{Reg} %
\DeclareMathOperator{\diff}{diff} %
\DeclareMathOperator{\bdim}{\mathbf{dim}} %

\newcommand{\ol}{\overline}

\newtheorem*{coro}{Corollary} 
\newtheorem*{lemm}{Lemma} 
\newtheorem*{prop}{Proposition} 
\newtheorem{theo}{Theorem}

\newtheorem{propn}{Proposition}[subsection] 

\def\ffrac#1#2{{\textstyle\frac{#1}{#2}}}




\title%
  [Geometry of homogeneous modules] 
  {Geometry and the zero sets of semi-invariants for homogeneous modules over canonical algebras}

\author{Grzegorz Bobi\'nski}

\address{Faculty of Mathematics and Computer Science \\ Nicolaus
Copernicus University \\ ul.~Chopina 12/18 \\ 87-100 Toru\'n \\
Poland} 

\email{gregbob@mat.uni.torun.pl} 

\begin{document} 

\begin{abstract}
We characterize the canonical algebras such that for all dimension
vectors of homogeneous modules the corresponding module varieties
are complete intersections (respectively, normal). We also
investigate the sets of common zeros of semi-invariants of
non-zero degree in important cases. In particular, we show that
for sufficiently big vectors they are complete intersections and
calculate the number of their irreducible components.
\end{abstract}



\maketitle

Throughout the paper $k$ denotes a fixed algebraically closed
field. By $\bbN$ and $\bbZ$ we denote the sets of nonnegative
integers and integers, respectively. Finally, if $i, j \in \bbZ$,
then $[i, j] = \{ l \in \bbZ \mid i \leq l \leq j \}$.

\section*{Introduction and main result}

Canonical algebras were introduced by Ringel
in~\cite{Rin1984}*{3.7} (for a definition see also~\ref{can_def}).
A canonical algebra $\Lambda$ depends on a sequence $(m_1, \ldots,
m_n)$, $n \geq 3$, of positive integers greater then $1$ and on a
sequence $(\lambda_3, \ldots, \lambda_n)$ of pairwise distinct
nonzero elements of $k$. In this situation we say that $\Lambda$
is a canonical algebra of type $(m_1, \ldots, m_n)$. The canonical
algebras play a prominent role in the representation theory of
algebras (see for example~\cites{LenPen1999, Sko1996}). In
particular, the module categories over canonical algebras are
derived equivalent to the categories of coherent sheaves over
weighted projective lines (see~\cite{GeigLen1987}). Moreover,
according to~\cite{Hap2001}*{Theorem~3.1} every quasi-titled
algebra is derived equivalent to a hereditary algebra or to a
canonical one.

Let $\Lambda$ be an algebra. For each element $\bd$ of the
Grothendieck group of $\Lambda$ one defines a variety
$\mod_\Lambda (\bd)$ called the variety of $\Lambda$-modules of
dimension vector $\bd$ (see~\ref{var_def}). The study of varieties
of modules is an important and interesting topic in the
representation theory of algebras (for some reviews of results see
for example~\cites{Bon1998, Geiss1996, Kra1982}).
In~\cite{BobSko1999} Skowro\'nski and the author proved that if
$\Lambda$ is a tame canonical algebra and $\bd$ is the dimension
vector of an indecomposable $\Lambda$-module, then $\mod_\Lambda
(\bd)$ is a complete intersection with at most~$2$ irreducible
components. The module varieties over canonical algebras were also
studied by Barot and Schr\"oer in~\cite{BarSch2001}.

Let $\Lambda$ be a canonical algebra of type $(m_1, \ldots, m_n)$.
We call a module regular if it is periodic with respect to the
action of the Auslander--Reiten translate (such modules are of
special interest in the representation theory, see for
example~\cite{Sko1999}). This class of modules also received
special attention from a geometric point of view. Skowro\'nski and
the author showed in~\cite{BobSko2002} that if $\bd$ is the
dimension vector of a regular module over a tame canonical algebra
$\Lambda$, then the corresponding variety is a normal complete
intersection. This result was extended in~\cite{Bob2005} by
showing that the varieties $\mod_\Lambda (\bd)$ are normal
(respectively, complete intersections) for all dimension vectors
$\bd$ of regular $\Lambda$-modules if and only if
\[
\ffrac{1}{m_1 - 1} + \cdots + \ffrac{1}{m_n - 1} > 2 n - 5 \quad
(\geq 2 n - 5).
\]

A special type of regular modules are the homogeneous ones, which
are invariant with respect to the action of the Auslander--Reiten
translate. The first result of the paper is the following.

\begin{theo} \label{theofirst}
Let $\Lambda$ be a canonical algebra of type $(m_1, \ldots, m_n)$.
\begin{enumerate}

\item
The varieties $\mod_\Lambda (\bd)$ are complete intersections for
all dimension vectors $\bd$ of homogeneous $\Lambda$-modules if
and only if
\[
\ffrac{1}{m_1} + \cdots + \ffrac{1}{m_n} \geq n - 4.
\]

\item
The varieties $\mod_\Lambda (\bd)$ are normal for all dimension
vectors $\bd$ of homogeneous $\Lambda$-modules if and only if
\[
\ffrac{1}{m_1} + \cdots + \ffrac{1}{m_n} > n - 4.
\]

\end{enumerate}
\end{theo}

If $\ffrac{1}{m_1} + \cdots + \ffrac{1}{m_n} \geq n - 2$, then the
algebra $\Lambda$ is tame, hence in this case the assertion
follows from the quoted result~\cite{BobSko2002}. Thus we may
assume that
\[
\ffrac{1}{m_1} + \cdots + \ffrac{1}{m_n} < n - 2.
\]
In this case the algebra $\Lambda$ is wild and $\bd$ is the
dimension vector of a homogeneous module if and only if $\bd = p
\bh$ for some $p \in \bbN$, where $\bh$ is the dimension vector
with all the coordinates equal to $1$.

We take a closer look into the boundary situation.

\begin{theo} \label{theosecond}
Let $\Lambda$ be a canonical algebra of type $(m_1, \ldots, m_n)$
with
\[
\ffrac{1}{m_1} + \cdots + \ffrac{1}{m_n} = n - 4,
\]
and let $m$ be the least common multiple of $m_1$, \ldots, $m_n$.
\begin{enumerate}

\item
If $m$ divides $p$, then $\mod_\Lambda (p \bh)$ is a complete
intersection with exactly two irreducible components.

\item
If $m$ does not divide $p$, then $\mod_\Lambda (p \bh)$ is a
normal complete intersection.

\end{enumerate}
\end{theo}

If $\Lambda$ is an algebra then for each dimension vector $\bd$ a
product $\GL (\bd)$ of general linear groups acts on the variety
$\mod_\Lambda (\bd)$ (see~\ref{var_def}). This action induces an
action on the ring $k [\mod_\Lambda (\bd)]$ of regular functions
on $\mod_\Lambda (\bd)$ (see~\ref{act_def}). It is known that for
a triangular algebra (no cycles in the Gabriel quiver), hence in
particular for canonical algebras, only the constant functions are
invariant with respect to this action, however the ring $\SI
[\mod_\Lambda (\bd)]$ of semi-invariants has a richer structure
(see for example~\cites{Hap1984, Rin1980, SkoWey2000}). In
particular, rings of semi-invariants arising for regular modules
over canonical algebras were studied~\cites{DomLen2000,
DomLen2002, SkoWey1999}.

In connection with rings of semi-invariants one may also ask, for
a dimension vector $\bd$ over an algebra $\Lambda$, about
properties of the set $\calZ (\bd)$ of the common zeros of the
semi-invariants of non-zero weight. This line of research (in the
context of module varieties) was initiated by Chang and Weyman
(\cite{ChaWey2004}) and continued by Riedtmann and Zwara
(\cites{RieZwa2003, RieZwa2004, RieZwa2006}). A motivation for
this research is that $\calZ (\bd)$ reflects properties of $k
[\mod_\Lambda (\bd)]$ as a module over $\SI [\mod_\Lambda (\bd)]$.

The last result of the paper concerns this topic.

\begin{theo} \label{theozeroset}
Let $\Lambda$ be a canonical algebra of type $(m_1, \ldots, m_n)$.
If
\[
\ffrac{1}{m_1} + \cdots + \ffrac{1}{m_n} < n - 4,
\]
then there exists $N$ such that, for $p \geq N$, $\calZ (p \bh)$
is a set theoretic complete intersection with
\[
(p - n) m_1 \cdots m_n + \sum_{l \in [1, n - 1]} \sum_{i_1 <
\cdots < i_l \in [1, n]} m_{i_1} \cdots m_{i_l} + 1
\]
irreducible components and $k [\mod_\Lambda (\bd)]$ is a free $\SI
[\mod_\Lambda (\bd)]$-module.
\end{theo}

Explicit bounds for $N$ can be found in
Propositions~\ref{prop_dom}, \ref{prop_tub}, and~\ref{prop_wild}.

The paper is organized as follows. In Section~\ref{sect_can} we
present definition and necessary facts about canonical algebras,
and in Section~\ref{sect_theo12} we prove Theorems~\ref{theofirst}
and~\ref{theosecond}. Next, in Section~\ref{sect_inv}, we collect
useful facts about semi-invariants, which in
Section~\ref{sect_theo3} are used in the proof of
Theorem~\ref{theozeroset}.

For basic background on the representation theory of algebras we
refer to~\cites{AssSimSko2006,AusReiSma1995}. Basic algebraic
geometry used in the article can be found for example
in~\cite{Kun}. The author gratefully acknowledges the support from
the Polish Scientific Grant KBN No.~1 P03A 018 27 and the
Schweizerischer Nationalfonds. The research leading to the results
presented in this paper was initiated while the author held a one
year post-doc position at the University of Bern.

\makeatletter
\def\@secnumfont{\mdseries} 
\makeatother

\section{Preliminaries on canonical algebras} \label{sect_can}

In this section we present facts about canonical algebras
necessary in our proofs.

\subsection{} 
\label{can_def} %
Let $\bm = (m_1, \ldots, m_n)$, $n \geq 3$, be a sequence of
integers greater than $1$ and let $\blambda = (\lambda_3, \ldots,
\lambda_n)$ be a sequence of pairwise distinct non-zero elements
of $k$. We define $\Lambda (\bm, \blambda)$ as the path algebra of
the bound quiver $(\Delta (\bm), R (\bm, \blambda))$, where
$\Delta (\bm)$ is the quiver
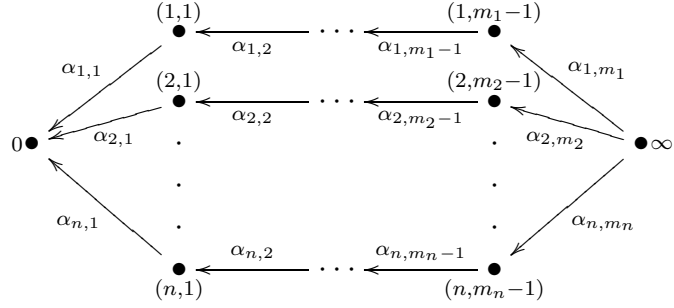
\[
\xymatrix@R=0.25\baselineskip@C=3\baselineskip{%
& \bullet \save*+!D{\scriptstyle (1, 1)} \restore
\ar[lddd]_{\alpha_{1, 1}} & \cdots \ar[l]^-{\alpha_{1, 2}} &
\bullet \save*+!D{\scriptstyle (1, m_1 - 1)} \restore
\ar[l]^-{\alpha_{1, m_1 - 1}}
\\ \\ %
& \bullet \save*+!D{\scriptstyle (2, 1)} \restore
\ar[ld]^{\alpha_{2, 1}} & \cdots \ar[l]^-{\alpha_{2, 2}} & \bullet
\save*+!D{\scriptstyle (2, m_2 - 1)} \restore \ar[l]^-{\alpha_{2,
m_2 - 1}}
\\ %
\bullet \save*+!R{\scriptstyle 0} \restore & \cdot & & \cdot &
\bullet \save*+!L{\scriptstyle \infty} \restore
\ar[luuu]_{\alpha_{1, m_1}} \ar[lu]^{\alpha_{2, m_2}}
\ar[lddd]^{\alpha_{n, m_n}}
\\ %
& \cdot & & \cdot
\\ %
& \cdot & & \cdot
\\ %
& \bullet \save*+!U{\scriptstyle (n, 1)} \restore
\ar[luuu]^{\alpha_{n, 1}} & \cdots \ar[l]_-{\alpha_{n, 2}} &
\bullet \save*+!U{\scriptstyle (n, m_n - 1)} \restore
\ar[l]_-{\alpha_{n, m_n - 1}} }
\]
and
\[
R (\bm, \blambda) = \{ \alpha_{1, 1} \cdots \alpha_{1, m_1} +
\lambda_i \alpha_{2, 1} \cdots \alpha_{2, m_2} - \alpha_{i, 1}
\cdots \alpha_{i, m_i} \mid i \in [3, n] \}.
\]
The algebras of the above form are called canonical. In
particular, we call $\Lambda (\bm, \blambda)$ a canonical algebra
of type $\bm$. If $\bm$ and $\blambda$ are fixed, then we usually
write $\Lambda$ and $(\Delta, R)$, instead of $\Lambda (\bm,
\blambda)$ and $(\Delta (\bm), R (\bm, \blambda))$, respectively.
In this case we further denote by $\Delta_0$ the set of vertices
of $\Delta$. Until the end of the section we assume that $\Lambda
= \Lambda (\bm, \blambda)$ is a fixed canonical algebra and
$(\Delta, R)$ is the corresponding bound quiver. The following
invariant
\[
\delta = \ffrac{1}{2} \bigl(n - 2 - \ffrac{1}{m_1} - \cdots -
\ffrac{1}{m_n} \bigr)
\]
controls the representation type of $\Lambda$
(see~\cite{Rin1984}). Namely, $\Lambda$ is tame if and only if
$\delta \leq 0$. Moreover, $\Lambda$ is domestic if and only if
$\delta < 0$.

\subsection{} 
By a representation of the bound quiver $(\Delta, R)$ we mean a
collection $M = (M_x, M_{i, j})_{x \in \Delta_0, \, i \in [1, n],
\, j \in [1, m_i]}$ of finite dimensional vector spaces $M_x$, $x
\in \Delta_0$, and linear maps $M_{i, j} : M_{(i, j)} \to M_{(i, j
- 1)}$, $i \in [1, n]$, $j \in [1, m_i]$, such that
\[
M_{1, 1} \cdots M_{1, m_1} + \lambda_i M_{2, 1} \cdots M_{2, m_2}
- M_{i, 1} \cdots M_{i, m_i} = 0, \, i \in [3, n],
\]
where $M_{(i, 0)} = M_0$ and $M_{(i, m_i)} = M_\infty$ for $i \in
[1, n]$. The category of representations of $(\Delta, R)$ is
equivalent to the category of $\Lambda$-modules, and we identify
$\Lambda$-modules and representations of $(\Delta, R)$. For a
representation $M$ we define its dimension vector $\bdim M \in
\bbN^{\Delta_0}$ by $(\bdim M)_x = \dim_k M_x$, $x \in \Delta_0$.

\subsection{} 
The Ringel bilinear form $\langle -, - \rangle : \bbZ^{\Delta_0}
\times \bbZ^{\Delta_0} \to \bbZ$ is defined
\begin{multline*}
\langle \bd', \bd'' \rangle = d_0' d_0'' + \sum_{i \in [1, n], j
\in [1, m_i - 1]} d_{i, j}' d_{i, j}'' + d_\infty' d_\infty''
\\ %
- \sum_{i \in [1, n], j \in [1, m_i]} d_{i, j}' d_{i, j - 1}'' +
(n - 2) d_\infty' d_0'',
\end{multline*}
where we use the convention that $d_{i, 0} = d_0$ and $d_{i, m_i}
= d_\infty$ for $\bd \in \bbZ^{\Delta_0}$ and $i \in [1, n]$, and
$d_{i, j} = d_{(i, j)}$ for $i \in [1, n]$ and $j \in [1, m_i -
1]$, plays an important role in describing the representation
theory of $\Lambda$. It is known (see~\cite{Bon1983}*{2.2}), that
if $M$ and $N$ are $\Lambda$-modules, then
\[
\langle \bdim M, \bdim N \rangle = [M, N] - [M, N]^1 + [M, N]^2,
\]
where following Bongartz~\cite{Bon1994} we write $[M, N] = \dim_k
\Hom_\Lambda (M, N)$, $[M, N]^1 = \dim_k \Ext_\Lambda^1 (M, N)$,
and $[M, N]^2 = \dim_k \Ext_\Lambda^2 (M, N)$.

\subsection{} 
Let $\bh$ be the dimension vector with all the coordinates equal
to $1$, and
\[
\be_{i, 0} = \bh - (\be_{i, 1} + \cdots + \be_{i, m_i - 1})
\]
for $i \in [1, n]$, where $(\be_x)_{x \in \bbZ^{\Delta_0}}$ is the
standard basis of $\bbZ^{\Delta_0}$ and $\be_{i, j} = \be_{(i,
j)}$ for $i \in [1, n]$ and $j \in [1, m_i - 1]$. One easily
checks that
\begin{align*}
\langle \bd, \bh \rangle & = d_0 - d_\infty = - \langle \bh, \bd
\rangle,
\\ %
\langle \bd, \be_{i, j} \rangle & = d_{i, j} - d_{i, j + 1}, \; i
\in [1, n], \; j \in [0, m_i - 1],
\\ %
\langle \be_{i, j}, \bd \rangle & = d_{i, j} - d_{i, j - 1}, \; i
\in [1, n], \; j \in [1, m_i - 1],
\\ %
\intertext{and}%
\langle \be_{i, 0}, \bd \rangle & = d_{i, m_i} - d_{i, m_i - 1},
\; i \in [1, n],
\end{align*}
for $\bd \in \bbZ^{\Delta_0}$.

\subsection{} 
\label{subcat}%
Let $\calP$ ($\calR$, $\calQ$, respectively) be the subcategory of
all $\Lambda$-modules which are direct sums of indecomposable
$\Lambda$-modules $X$ such that
\[
\langle \bdim X, \bh \rangle > 0 \quad (\langle \bdim X, \bh
\rangle = 0, \, \langle \bdim X, \bh \rangle < 0, \text{
respectively}).
\]
We have the following properties of the above decomposition of the
category of $\Lambda$-modules (see~\cite{Rin1984}*{3.7}).

First, $[N, M] = 0$ and $[M, N]^1 = 0$, if either $N \in \calR
\vee \calQ$ and $M \in \calP$, or $N \in \calQ$ and $M \in \calP
\vee \calR$. Here, for two subcategories $\calX$ and $\calY$ of
the category of $\Lambda$-modules, we denote by $\calX \vee \calY$
the additive closure of their union. Moreover, one knows that
$\pd_\Lambda M \leq 1$ for $M \in \calP \vee \calR$ and
$\id_\Lambda N \leq 1$ for $N \in \calR \vee \calQ$. Secondly,
$\calR$ decomposes into a $\bbP^1 (k)$-family $\coprod_{\lambda
\in \bbP^1 (k)} \calR_\lambda$ of uniserial categories. In
particular, $[M, N] = 0$ and $[M, N]^1 = 0$ if $M \in
\calR_\lambda$ and $N \in \calR_\mu$ for $\lambda \neq \mu$. We
put $\calR' = \coprod_{i \in [1, n]} \calR_{\lambda_i}$ and
$\calR'' = \coprod_{\lambda \in \bbP^1 (k) \setminus \{ \lambda_1,
\ldots, \lambda_n \}} \calR_\lambda$, where $\lambda_1 = 0$ and
$\lambda_2 = \infty$. If $\lambda \in \bbP^1 (k) \setminus \{
\lambda_1, \ldots, \lambda_n \}$, then there is a unique simple
object $R_\lambda$ in $\calR_\lambda$ and its dimension vector is
$\bh$. On the other hand, if $\lambda = \lambda_i$ for $i \in [1,
n]$, then there are $m_i$ simple objects $R_\lambda^{(0)}$,
\ldots, $R_\lambda^{(m_i - 1)}$ in $\calR_{\lambda_i}$ and their
dimension vectors are $\be_{i, j}$, $j \in [0, m_i - 1]$,
respectively.

\subsection{} 
\label{subsectP}%
Let $\bP$, $\bR$ and $\bQ$ denote the sets of the dimension vector
of modules from $\calP$, $\calR$, $\calQ$, respectively. According
to~\cite{Bob2005}*{2.6}, $\bd \in \bP$ if and only if either $\bd
= 0$ or $d_0 > d_\infty \geq 0$ and $d_{i, j} \geq d_{i, j + 1}$
for all $i \in [1, n]$ and $j \in [0, m_i - 1]$. Dually, $\bd \in
\bQ$ if and only if either $\bd = 0$ or $0 \leq d_0 < d_\infty$
and $d_{i, j} \leq d_{i, j + 1}$ for all $i \in [1, n]$ and $j \in
[0, m_i - 1]$.

For $l_1 \in [0, m_1 - 1]$, \ldots, $l_n \in [0, m_n - 1]$ we set
\[
\be (l_1, \ldots, l_n) = \be_0 + \sum_{i \in [1, n]} \sum_{j \in
[1, l_i]} \be_{i, j}.
\]
We will need the following fact.

\begin{lemm} \label{lemmdP}
If $\bd \in \bP$ is such that $\langle \bd, \bh \rangle = 1$, then
\[
\bd = r \bh + \be (l_1, \ldots, l_n)
\]
for some $r \in \bbN$ and $l_i \in [0, m_i - 1]$, $i \in [1, n]$.
In particular, $\langle \bd, \bd \rangle = 1$.
\end{lemm}

\begin{proof}
The first part follows immediately from the above description of
$\bP$ and the equality $\langle \bd, \bh \rangle = d_0 -
d_\infty$. The second part follows by direct calculations.
\end{proof}

\subsection{} 
The following inequality will be extremely useful in our proofs.

\begin{lemm} \label{lemmineq}
Let $\bd \in \bbZ^{\Delta_0}$. Then
\[
\langle \bd, \bd \rangle \geq -\delta (d_0 - d_\infty)^2.
\]
Moreover, the equality holds if and only if
\[
d_{i, j} = \ffrac{1}{m_i} ((m_i - j) d_0 + j d_\infty)
\]
for $i \in [1, n]$ and $j \in [1, m_i - 1]$.
\end{lemm}

\begin{proof}
Let $\bd' = \bd - d_\infty \bh$. Then $\langle \bd, \bd \rangle =
\langle \bd', \bd' \rangle$ and $d_0' = d_0 - d_\infty$, hence the
claim follows from the following equality
\begin{multline*}
\langle \bd', \bd' \rangle = - \delta d_0'^2 +
\\ %
\ffrac{1}{2} \sum_{i \in [1, n]} \sum_{j \in [1, m_i - 1]}
\ffrac{1}{(m_i - j) (m_i - j + 1)} ((m_i - j + 1) d_{i, j}' - (m_i
- j) d_{i, j - 1}')^2,
\end{multline*}
which was suggested to me by Professor Riedtmann.
\end{proof}

\makeatletter
\def\@secnumfont{\mdseries} 
\makeatother

\section{Preliminaries on module varieties
\\ 
and proofs of Theorems~\ref{theofirst} and~\ref{theosecond}}
\label{sect_theo12}

Throughout this section $\Lambda$ is a fixed canonical algebra of
type $\bm$ and $\Delta$ is its quiver. We first define
in~\ref{var_def} varieties of modules, then formulate
in~\ref{crit_geom} numerical criteria characterizing geometric
properties of these varieties, and finally we apply these criteria
in~\ref{proof_one} and~\ref{proof_two} in order to prove
Theorems~\ref{theofirst} and~\ref{theosecond}, respectively.

\subsection{} 
\label{var_def}%
For $\bd \in \bbN^{\Delta_0}$ let $\bbA (\bd) = \prod_{i \in [1,
n], j \in [1, m_i]} \bbM (d_{i, j - 1}, d_{i, j})$. By
$\mod_\Lambda (\bd)$ we denote the subset of $\bbA (\bd)$ formed
by all tuples $(M_{i, j})$ such that
\[
M_{1, 1} \cdots M_{1, m_1} + \lambda_i M_{2, 1} \cdots M_{2, m_2}
- M_{i, 1} \cdots M_{i, m_i} = 0, \, i \in [3, n].
\]
We identify the points $M$ of $\mod_\Lambda (\bd)$ with
$\Lambda$-modules of dimension vector $\bd$ by taking $M_x =
k^{d_x}$ for $x \in \Delta_0$. The product $\GL (\bd) = \prod_{x
\in \Delta_0} \GL (d_x)$ of general linear groups acts on
$\mod_\Lambda (\bd)$ by conjugation
\[
(g \cdot M)_{i, j} = g_{(i, j - 1)} M_{i, j} g_{(i, j)}^{-1}, \; i
\in [1, n], \, j \in [1, m_i],
\]
for $g \in \GL (\bd)$ and $M \in \mod_\Lambda (\bd)$, where
$g_{(i, 0)} = g_0$ and $g_{(i, m_i)} = g_\infty$ for $i \in [1,
n]$. The orbits with respect to this action correspond bijectively
to the isomorphism classes of $\Lambda$-modules of dimension
vector $\bd$. For $M \in \mod_\Lambda (\bd)$ we denote by $\calO
(M)$ the $\GL (\bd)$-orbit of $M$. We put
\[
a (\bd) = \dim \bbA (\bd) - (n - 2) d_0 d_\infty.
\]
Note that $a (\bd) = \dim \GL (\bd) - \langle \bd, \bd \rangle$.

\subsection{} 
\label{crit_geom}%
For a subcategory $\calX$ of the category of $\Lambda$-modules and
a dimension vector $\bd$ we denote by $\calX (\bd)$ the set of all
$M \in \mod_\Lambda (\bd)$ such that $M \in \calX$. One knows that
if $\bd \in \bbN^{\Delta_0}$ then $\calP (\bd)$ and $(\calR \vee
\calQ) (\bd)$ are open subsets of $\mod_\Lambda (\bd)$
(see~\cite{Bob2005}*{Lemmas~3.7 and~3.8}). Together with the
properties of the categories $\calP$, $\calR$ and $\calQ$ listed
in~\ref{subcat}, it implies that we can apply the results
of~\cite{Bob2005}*{Section~4} with $\calX = \calP$ and $\calY =
\calR \vee \calQ$.

Observe that if $\bd \in \bP$ and $p \geq d_0$, then $p \bh - \bd
\in \bQ$. Moreover,
\[
\langle p \bh - \bd, \bd \rangle = - p (d_0 - d_\infty) - \langle
\bd, \bd \rangle.
\]
Thus, the following proposition is a consequence of
\cite{Bob2005}*{Proposition~4.3, Proposition~4.5
and~Proposition~4.9}.

\begin{propn} \label{propfirst}
Let $p \geq 1$.
\begin{enumerate}

\item \label{pointoneone}
The variety $\mod_\Lambda (p \bh)$ is a complete intersection if
and only if $\langle \bd, \bd \rangle \geq - p (d_0 - d_\infty)$
for all $\bd \in \bP$ such that $d_0 \leq p$.

\item
The variety $\mod_\Lambda (p \bh)$ is normal if and only if
$\langle \bd, \bd \rangle > - p (d_0 - d_\infty)$ for all $\bd \in
\bP$, $\bd \neq 0$, such that $d_0 \leq p$.

\end{enumerate}
\end{propn}

We will also need the following consequence of the proof
of~\cite{Bob2005}*{Proposition~4.5}.

\begin{propn} \label{propirr}
Let $p \geq 1$ and assume that $\langle \bd, \bd \rangle \geq - p
(d_0 - d_\infty)$ for all $\bd \in \bP$ such that $d_0 \leq p$.
Then the irreducible components of $\mod_\Lambda (p \bh)$ are in
bijection with the dimensions vectors $\bd \in \bP$ such that $d_0
\leq p$ and $\langle \bd, \bd \rangle = - p (d_0 - d_\infty)$.
\end{propn}

\subsection{} 
\label{proof_one}%
We now prove Theorem~\ref{theofirst}. Recall that it is enough to
prove the theorem for $\bd = p \bh$ with $p > 0$ (see the
discussion after Theorem~\ref{theofirst}). Assume that
\[
\ffrac{1}{m_1} + \cdots + \ffrac{1}{m_n} \geq n - 4.
\]
Then $\delta \leq 1$ and according to Lemma~\ref{lemmineq}
\[
\langle \bd, \bd \rangle \geq - \delta (d_0 - d_\infty)^2 \geq -
d_0 (d_0 - d_\infty) \geq - p (d_0 - d_\infty)
\]
for each $p \geq 1$ and $\bd \in \bP$ such that $d_0 \leq p$.
According to Proposition~\ref{propfirst}\eqref{pointoneone}, this
implies that $\mod_\Lambda (p \bh)$ is a complete intersection.
Analogously, we prove that $\mod_\Lambda (p \bh)$ is normal if
\[
\ffrac{1}{m_1} + \cdots + \ffrac{1}{m_n} > n - 4,
\]
since in this case $\delta < 1$ and the second inequality in the
above string is strict for $\bd \neq 0$.

It remains to prove that if
\[
\ffrac{1}{m_1} + \cdots + \ffrac{1}{m_n} < n - 4 \qquad (\leq n -
4),
\]
then there exists $p$ such that $\mod_\Lambda (p \bh)$ is not a
complete intersection (respectively, normal), or in other words,
there exists $\bd \in \bP$ ($\bd \neq 0$) such that
\[
\langle \bd, \bd \rangle < - p (d_0 - d_\infty) \qquad (\leq p
(d_0 - d_\infty)).
\]
A construction of such $p$ and $\bd$ is suggested by
Lemma~\ref{lemmineq}. Namely, let $p = m_1 \cdots m_n$ and $\bd$
be given by the formulas
\begin{align*}
d_0 & = p, & d_\infty & = 0, & d_{i, j} & = \ffrac{m_i - j}{m_i}
p, \; i \in [1, n], \, j \in [1, m_i - 1].
\end{align*}
Then $\bd \in \bP$, $\bd \neq 0$, and $\langle \bd, \bd \rangle =
- \delta p (d_0 - d_\infty)$, what finishes the proof.

\subsection{} 
\label{proof_two}%
We now prove Theorem~\ref{theosecond}. Assume that
\[
\ffrac{1}{m_1} + \cdots + \ffrac{1}{m_n} = n - 4
\]
and $p > 0$. We already know that $\mod_\Lambda (p \bh)$ is a
complete intersection. Moreover, according
to~\cite{Bob2005}*{Proposition ~4.9} it is normal if and only it
is irreducible. Thus our task it to classify the irreducible
components of $\mod_\Lambda (p \bh)$. According to
Proposition~\ref{propirr} this is equivalent to classifying the
dimension vectors $\bd \in \bP$ such that $d_0 \leq p$ and
$\langle \bd, \bd \rangle = - p (d_0 - d_\infty)$. Obviously, one
such vector is the zero vector. Hence assume that $\bd \neq 0$. It
follows from Lemma~\ref{lemmineq} that $\langle \bd, \bd \rangle
\geq - (d_0 - d_\infty)^2$ (recall that $\delta = 1$ in our case).
Thus the condition $\langle \bd, \bd \rangle = - p (d_0 -
d_\infty)$ implies that $d_0 = p$ and $d_\infty = 0$. Moreover, we
know again from Lemma~\ref{lemmineq} that $\langle \bd, \bd
\rangle = -(d_0 - d_\infty)^2$ if and only if $d_{i, j} =
\ffrac{m_i - i }{m_i} p$, $i \in [1, n]$, $j \in [1, m_i - 1]$.
Note that $\bd$ defined by the above formulas belongs to $\bP$ if
and only if it belongs to $\bbN^{\Delta_0}$, i.e., if and only if
$m_i$ divides $p$ for $i \in [1, n]$. This observation concludes
the proof.

\makeatletter
\def\@secnumfont{\mdseries} 
\makeatother

\section{Preliminaries on semi-invariants} \label{sect_inv}

Throughout this section $\Lambda = \Lambda (\bm, \blambda)$ is a
fixed canonical algebra, and $\Delta = \Delta (\bm)$. Moreover, we
put
\[
|\bm| = m_1 + \cdots + m_n.
\]
Our main aim in this section is to prove
Proposition~\ref{propcompbis}, which reduces the proof of
Theorem~\ref{theozeroset} to a certain inequality. In order to
achieve this aim we first recall basic facts about semi-invariants
in~\ref{act_def} and~\ref{construction}. The main result of this
first part is Corollary~\ref{coroZ} giving a new formulation of
Theorem~\ref{theozeroset}, which we subsequently improve
in~\ref{lemmRbis}--\ref{propcompbis}.

\subsection{} 
\label{act_def}%
Let $\bd \in \bbN^{\Delta_0}$. The action of $\GL (\bd)$ on
$\mod_\Lambda (\bd)$ induces an action of $\GL (\bd)$ on the
coordinate ring $k [\mod_\Lambda (\bd)]$ of $\mod_\Lambda (\bd)$
in the usual way, i.e.\
\[
(g \cdot f) (M) = f (g^{-1} \cdot M)
\]
for $g \in \GL (\bd)$, $f \in k [\mod_\Lambda (\bd)]$, and $M \in
\mod_\Lambda (\bd)$. The product $\SL (\bd) = \prod_{x \in
\Delta_0} \SL (d_x)$ of special linear groups is a closed subgroup
of $\GL (\bd)$. The ring $\SI [\mod_\Lambda (\bd)]$ of invariants
with respect to the induced action of $\SL (\bd)$ on $k
[\mod_\Lambda (\bd)]$ is called the ring of semi-invariants.

By a weight we mean a group homomorphism $\sigma : \bbZ^{\Delta_0}
\to \bbZ$. We identify the weights with the elements of the group
$\bbZ^{\Delta_0}$ in the usual way. If $\sigma$ is a weight, then
we define the weight space
\[
\SI [\mod_\Lambda (\bd)]_\sigma = \{ f \in k [\mod_\Lambda (\bd)]
\mid g \cdot f = \bigl( \prod_{x \in \Delta_0}
\det\nolimits^{\sigma (x)} (g) \bigr) f \}
\]
(observe that $\SI [\mod_\Lambda (\bd)]_\sigma \subset \SI
[\mod_\Lambda (\bd)]$). It is known that
\[
\SI [\mod_\Lambda (\bd)] = \bigoplus_{\sigma \in \bbZ^{\Delta_0}}
\SI [\mod_\Lambda (\bd)]_\sigma
\]
provided $\bd$ is sincere, i.e.\ $d_x \neq 0$ for $x \in
\Delta_0$. Moreover $\SI [\mod_\Lambda (\bd)]_0 = k$. In this
situation the set $\calZ (\bd)$ of common zeros of homogeneous
semi-invariants with non-zero weights is called the zero set of
semi-invariants.

\subsection{} 
\label{construction}%
We recall now a construction of semi-invariants described
in~\cite{DomLen2002} (being a generalization of a construction of
Schofield~\cite{Sch1991} --- compare also~\cites{Dom2002,
DerWey2002}). Let $M$ be a $\Lambda$-module of projective
dimension at most $1$, and let
\[
0 \to P_1 \xrightarrow{\varphi} P_0 \to M \to 0
\]
be its minimal projective resolution. If $\bd \in \bbN^{\Delta_0}$
satisfies $\langle \bdim M, \bd \rangle = 0$, then the map
$d_{\bd}^M : \mod_\Lambda (\bd) \to k$ given by $d_{\bd}^M (N) =
\det \Hom_\Lambda (\varphi, N)$ is well-defined (up to scalars)
and is a homogeneous semi-invariant of weight $-\langle \bdim M, -
\rangle$. Moreover, $d_{\bd}^M (N) = 0$ if and only if $[M, N]
\neq 0$.

\subsection{} 
For $\bd \in \bbN^{\Delta_0}$ let $\Reg_\Lambda (\bd)$ denote the
closure of $\calR (\bd)$. If $\calR (\bd)$ is nonempty, i.e.\ $\bd
\in \bR$, then $\Reg_\Lambda (\bd)$ is an irreducible component of
$\mod_\Lambda (\bd)$. The action of $\GL (\bd)$ on $k
[\mod_\Lambda (\bd)]$ restricts to an action on $k [\Reg_\Lambda
(\bd)]$. In particular, by $\SI [\Reg_\Lambda (\bd)]$ we denote
the ring of $\SI (\bd)$-invariant regular functions on
$\Reg_\Lambda (\bd)$. The rings $\SI [\Reg_\Lambda (\bd)]$ for
$\bd \in \bR$ have been studied in~\cite{SkoWey1999} (in case of
characteristic $0$) and in~\cite{DomLen2002} (in case of arbitrary
characteristic). We now list their properties which are important
for our investigations.

\begin{prop}
Let $p \geq 1$. Then $\SI [\Reg_\Lambda (p \bh)]$ is generated by
$d_{p \bh}^M$, $M \in \calR$. Moreover, if $p \geq n - 1$ then
$\SI [\Reg_\Lambda (p \bh)]$ is a polynomial algebra in $|\bm| + p
+ 1 - n$ variables.
\qed 
\end{prop}

The following consequence of the above proposition will be crucial
for us.

\begin{coro} \label{coroZ}
Let $p \geq 1$. If $\mod_\Lambda (p \bh)$ is irreducible, then
\[
\calZ (p \bh) = \{ N \in \mod_\Lambda (p \bh) \mid [M, N] \neq 0
\text{ for all } M \in \calR \}.
\]
Moreover, if $p \geq n - 1$ and
\[
\dim \calZ (p \bh) = a (\bd) - |\bm| - p - 1 + n,
\]
then $\calZ (p \bh)$ is a set theoretic complete intersection and
$k [\mod_\Lambda (p \bh)]$ is a free $\SI [\mod_\Lambda (p
\bh)]$-module.
\end{coro}

\begin{proof}
Recall from~\cite{Bob2005} that if $\mod_\Lambda (p \bh)$ is
irreducible then $\mod_\Lambda (p \bh)$ is a complete intersection
of dimension $a (\bd)$, thus the above corollary is a direct
consequence of the above proposition
and~\cite{Bob2006}*{Section~4}.
\end{proof}

We note that always $\calZ (p \bh) \geq a (\bd) - |\bm| - p - 1 +
n$, hence the hard part of the proof is show that $\calZ (p \bh)
\leq a (\bd) - |\bm| - p - 1 + n$.

\subsection{} 
We will need the following well-known fact.

\begin{lemm} \label{lemmRbis}
If $p \geq 1$, then $\calR'' (p \bh)$ is an open subset of $\calR
(p \bh)$. In particular, $\dim \calR'' (p \bh) = a (p \bh)$.
\end{lemm}

\begin{proof}
Let $M \in \calR$. Then $M \in \calR''$ if and only if
$\Hom_\Lambda (R_{\lambda_i}^{(j)}, M) = 0$ for all $i \in [1, n]$
and $j \in [0, m_i - 1]$, which implies the first part of the
lemma. The second part is an immediate consequence of the
well-known fact that $\calR (p \bh)$ is an irreducible set of
dimension $a (p \bh)$ (see for example remarks after
\cite{Bob2005}*{Lemma~3.7}).
\end{proof}

\subsection{} 
Fix $\bd' \in \bP$, $\bd'' \in \bQ$ and $X \in \calR'$ such that
$\bd' + \bd'' + \bdim X = q \bh$ for some $q \geq 1$. For each $p
\geq q$ we consider the set $\calC_p (\bd', \bd'', X)$ consisting
of all $M \in \mod_\Lambda (p \bh)$ which are isomorphic to
modules of the form $M' \oplus M'' \oplus X \oplus Y$, where $M'
\in \calP (\bd')$, $M'' \in \calQ (\bd'')$ and $Y \in \calR'' ((p
- q) \bh)$. We will need the following properties of the set
$\calC_p (\bd', \bd'', X)$.

\begin{lemm}
Let $\bd'$, $\bd''$, $X$, $q$ and $p$ be as above, and $\calC =
\calC_p (\bd', \bd'', X)$. Then $\calC$ is an irreducible
constructible set of dimension
\[
a (p \bh) - ((2 p - q) \langle \bd', \bh \rangle + \langle \bd',
\bd' \rangle + \langle \bd', \bdim X \rangle + [X, X]).
\]
\end{lemm}

\begin{proof}
The claim follows from \cite{Bob2005}*{Corollary~3.4}. Indeed
\[
\calC = \calP (\bd') \oplus \calO (X) \oplus \calR'' ((p - q) \bh)
\oplus \calQ (\bd'')
\]
in the notation of \cite{Bob2005}*{3.4}. Moreover, according to
\cite{Bob2005}*{Lemma~3.8}
\[
\dim \calP (\bd') = a (\bd') \qquad \text{ and } \qquad \dim \calQ
(\bd'') = a (\bd'').
\]
Further, by a well-known formula for the dimension of $\calO (X)$
(see for example~\cite{KraRie}*{2.2})
\[
\dim \calO (X) = \dim \GL (\bd) - [X, X] = a (\bd) + \langle \bd,
\bd \rangle - [X, X],
\]
where $\bd = \bdim X$. In addition, according to
Lemma~\ref{lemmRbis}
\[
\dim \calR'' ((p - q) \bh) = a ((p - q) \bh).
\]
Finally, for any $M' \in \calP (\bd')$, $M'' \in \calQ (\bd'')$
and $M \in \calR'' ((p - q) \bh)$
\begin{align*}
[M', X] & = \langle \bd', \bd \rangle, & [X, M'] & = 0,
\\ %
[M', M] & = \langle \bd', (p - q) \bh \rangle, & [M, M'] & = 0,
\\ %
[M', M''] & = \langle \bd', \bd'' \rangle, & [M', M''] & = 0,
\\ %
[X, M] & = 0 = \langle \bd, (p - q) \bh \rangle, & [M, X] & = 0,
\\ %
[X, M''] & = \langle \bd, \bd'' \rangle, & [M'', X] & = 0,
\\ %
[M, M''] & = \langle (p - q) \bh, \bd'' \rangle, & [M'', M] & = 0,
\end{align*}
hence we are in position to apply~\cite{Bob2005}*{Corollary~3.4}.
Since $\langle \bh, \bd' \rangle = - \langle \bd', \bh \rangle$
and $\langle \bd'', \bh \rangle = \langle q \bh - \bd - \bd', \bh
\rangle = - \langle \bd', \bh \rangle$, the formula follows by
direct calculations.
\end{proof}

\subsection{} 
Another important property is the following.

\begin{lemm}
Let $\bd'$, $\bd''$, $X$, $q$ and $p$ be as above, and $\calC =
\calC_p (\bd', \bd'', X)$. If $\mod_\Lambda (p \bh)$ is
irreducible, then $\calC \cap \calZ (p \bh) \neq \varnothing$ if
and only if $\bd' \neq 0$ and for each $i \in [1, n]$ and $j \in
[0, m_i - 1]$ either $\langle \bd', \be_{i, j} \rangle \neq 0$ or
$[X, R_{\lambda_i}^{(j)}] \neq 0$. In particular, $\calC \cap
\calZ (p \bh) \neq \varnothing$ if and only if $\calC \subset
\calZ (p \bh)$.
\end{lemm}

\begin{proof}
Take $N \in \calC$. Recall that, according to
Corollary~\ref{coroZ}, under the assumptions of the lemma $N \in
\calZ (p \bh)$ if and only if $[N, R] \neq 0$ for all $R \in
\calR$. This is equivalent to saying that $[N, R_\lambda] \neq 0$
for all $\lambda \in \bbP^1 (k) \setminus \{ \lambda_1, \ldots,
\lambda_n \}$ and $[N, R_{\lambda_i}^{(j)}] \neq 0$ for all $i \in
[1, n]$ and $j \in [0, m_i - 1]$. Write $N \simeq M' \oplus X
\oplus M \oplus M''$ for $M' \in \calP$, $M \in \calR''$ and $M
\in \calQ$. The former condition is equivalent to $M' \neq 0$
(i.e., $\bd' \neq 0$) since $[X \oplus M'', R] = 0$ for all $R \in
\calR''$ and $[M, R_\lambda] = 0$ for all but a finite number of
$\lambda \in \bbP^1 (k) \setminus \{ \lambda_1, \ldots, \lambda_n
\}$. Similarly, the latter condition leads to the second condition
of the lemma.
\end{proof}

\subsection{} 
For $p \geq 1$ let $\fZ_p$ be the set of all triples $(\bd',
\bd'', [X])$, where $\bd' \in \bP$, $\bd'' \in \bQ$ and $X \in
\calR'$ are such that $\bd' + \bd'' + \bdim X = q \bh$ for some $q
\leq p$, $\bd' \neq 0$, and for each $i \in [1, n]$ and $j \in [0,
m_i - 1]$ either $\langle \bd', \be_{i, j} \rangle \neq 0$ or $[X,
R_{\lambda_i}^{(j)}] \neq 0$. Here $[X]$ denotes the isomorphism
class of $X$. Observe that $\fZ_p$ is a finite set, hence as a
consequence of two preceding lemmas and Corollary~\ref{coroZ} we
get the following.

\begin{prop} \label{propcomp}
Let $p \geq 1$ and assume that $\mod_\Lambda (p \bh)$ is
irreducible. If
\renewcommand{\theequation}{$*$}
\begin{equation} \label{eqineq}
(2 p - q) \langle \bd', \bh \rangle + \langle \bd', \bd' \rangle +
\langle \bd', \bdim X \rangle + [X, X] \geq |\bm| + p + 1 - n
\end{equation}
for all $(\bd', \bd'', [X]) \in \fZ_p$, where $q \bh = \bd' +
\bd'' + \bdim X$, then $\calZ (p \bh)$ is a set theoretic complete
intersection  and $k [\mod_\Lambda (p \bh)]$ is free as a module
over $\SI [\mod_\Lambda (p \bh)]$. Moreover, if this is the case
then the map
\[
\fZ_p \ni (\bd', \bd'', [X]) \mapsto \ol{\calC_p (\bd', \bd'', X)}
\subset \calZ (p \bh)
\]
induces a bijection between those members of $\fZ_p$ with equality
in~\eqref{eqineq} and the irreducible components of $\calZ (p
\bh)$.
\end{prop}

\begin{proof}
The only missing part is the well-known fact that irreducible
components of complete intersections have the same dimension
(\cite{Kun}*{3.12}).
\end{proof}

\subsection{} 
The following inequality will give us a more accessible version of
the previous fact.

\begin{lemm} \label{lemmEndX}
Let $p \geq 1$. If $(\bd', \bd'', [X]) \in \fZ_p$, then
\[
[X, X] \geq |\bm| - n \langle \bd', \bh \rangle.
\]
\end{lemm}

\begin{proof}
Since the categories $\calR_{\lambda_i}$, $i \in [1, n]$, are
uniserial and pairwise orthogonal it follows for indecomposable $R
\in \calR'$, that if $[R, R_{\lambda_{i_0}}^{(j_0)}] \neq 0$ for
some $i \in [1, n]$ and $j \in [0, m_i - 1]$, then $[R,
R_{\lambda_i}^{(j)}] = 0$ for all $i \in [1, n]$ and $j \in [0,
m_i - 1]$ such that $(i, j) \neq (i_0, j_0)$. For $i \in [1, n]$
let $s_i$ denotes the number of the indecomposable direct summands
of $X$ which belong to $\calR_{\lambda_i}$. Since
\[
\langle \bd', \bh \rangle = \langle \bd', \be_{i, 0} \rangle +
\cdots + \langle \bd', \be_{i, m_i - 1} \rangle,
\]
for each $i \in [1, n]$, $\langle \bd', \be_{i, j} \rangle \geq 0$
for all $i \in [1, n]$ and $j \in [0, m_i - 1]$, it follows from
the definition of $\fZ_p$ that $\langle \bd', \bh \rangle \geq m_i
- s_i$. Using that $[X, X] \geq s_1 + \cdots + s_n$, we obtain our
claim.
\end{proof}

\subsection{} 
We now reformulate Proposition~\ref{propcomp}.

\begin{prop} \label{propcompbis}
Let $p \geq 1$ and assume that $\mod_\Lambda (p \bh)$ is
irreducible. If
\[
(p - q) \langle \bd', \bh \rangle + (p - n) (\langle \bd', \bh
\rangle - 1) + (\langle \bd', \bd' \rangle - 1) \geq 0
\]
for all $(\bd', \bd'', [X]) \in \fZ_p$, where $q \bh = \bd' +
\bd'' + \bdim X$, then $\calZ (p \bh)$ is a set theoretic complete
intersection and $k [\mod_\Lambda (p \bh)]$ is free as a module
over $\SI [\mod_\Lambda (p \bh)]$. Moreover, if this is the case
and the above inequality is strict for all $(\bd', \bd'', [X]) \in
\fZ_p$ such that $\langle \bd', \bh \rangle > 1$, then the
irreducible components are indexed by the triples $(\bd', \bd'',
[X]) \in
\fZ_p$ such that %
\renewcommand{\theequation}{$+$}
\begin{equation} \label{eqcond}
\begin{split}
\langle \bd', \bdim X \rangle = 0, &  \qquad [X, X] = |\bm| - n
\langle \bd', \bh \rangle,
\\ %
\langle \bd', \bh \rangle = 1, &  \qquad \bd' + \bd'' + \bdim X =
p \bh.
\end{split}
\end{equation}
\end{prop}

\begin{proof}
The first part follows from Proposition~\ref{propcomp} and
Lemma~\ref{lemmEndX} together with the obvious inequality $\langle
\bd', \bdim X \rangle \geq 0$. The second part is obtained in a
similar way: one has to use in addition Lemma~\ref{lemmdP}.
\end{proof}

For future reference we introduce the following notation:
\[
\diff (\bd', \bd'', [X]) = (p - q) \langle \bd', \bh \rangle + (p
- n) (\langle \bd', \bh \rangle - 1) + (\langle \bd', \bd' \rangle
- 1)
\]
for $(\bd', \bd'', [X]) \in \fZ_p$, with $\bd' + \bd'' + \bdim X =
q \bh$.

\subsection{} 
We calculate now the number of triples described in the above
proposition.

\begin{lemm} \label{lemmirr}
If $p \geq n$ then the number of triples $(\bd', \bd'', [X]) \in
\fZ_p$ satisfying~\eqref{eqcond} is
\[
(p - n) m_1 \cdots m_n + \sum_{l \in [1, n - 1]} \sum_{i_1 <
\cdots < i_l \in [1, n]} m_{i_1} \cdots m_{i_l} + 1.
\]
\end{lemm}

\begin{proof}
It follows from Lemma~\ref{lemmdP} that the condition $\langle
\bd', \bh \rangle = 1$ implies that
\[
\bd' = r \bh + \be_0 + \be (l_1, \ldots, l_n)
\]
for some $r \in \bbN$ and $l_i \in [0, m_i - 1]$, $i \in [1, n]$.
Note that if $i \in [1, n]$ and $j \in [0, m_i - 1]$, then
$\langle \bd', \be_{i, j} \rangle > 0$ if and only if $j = l_i$.
Consequently, for each $i \in [1, n]$ and $j \in [0, m_i - 1]$, $j
\neq l_i$, there exists an indecomposable direct summand $X_{i,
j}$ of $X$ such that $[X_{i, j}, R_{\lambda_i}^{(j)}] \neq 0$. It
follows that $X_{i, j} = R_{\lambda_i}^{(j)}$ since otherwise
either $[X_{i, l}, X_{i, j}] \neq 0$ for $l \neq j$ and
consequently $[X, X] > |\bm| - n$, or $\langle \bd', \bdim X_{i,
j} \rangle \neq 0$. It is possible to find $\bd'' \in \bQ$ such
that $\bd' + \bd'' + \bdim X = p \bh$ if and only if $r + |\{ i
\in [1, n] \mid l_i > 0 \}| \leq p - 1$, which implies the formula
in the lemma.
\end{proof}

\makeatletter
\def\@secnumfont{\mdseries} 
\makeatother

\section{Proof of Theorem~\ref{theozeroset}} \label{sect_theo3}
Throughout this section we assume that $\Lambda$ is a fixed
canonical algebra of type $\bm$. Our aim in this section is to
show how Proposition~\ref{propcompbis} and Lemma~\ref{lemmirr}
imply Theorem~\ref{theozeroset}.

\subsection{} 
We start with the domestic case.

\begin{prop} \label{prop_dom}
Let $\delta < 0$. If $p \geq n$ then $\calZ (p \bh)$ is a set
theoretic complete intersection and $k [\mod_\Lambda (\bd)]$ is a
free $\SI [\mod_\Lambda (\bd)]$-module. For $p > n$ the number of
the irreducible components of $\calZ (p \bh)$ is
\[
(p - n) m_1 \cdots m_n + \sum_{l \in [1, n - 1]} \sum_{i_1 <
\cdots < i_l \in [1, n]} m_{i_1} \cdots m_{i_l} + 1.
\]
\end{prop}

\begin{proof}
The claim follows from Proposition~\ref{propcompbis},
Lemma~\ref{lemmirr} and Theorem~\ref{theofirst}. It is enough to
observe that, according to Lemma~\ref{lemmineq}, $\langle \bd',
\bd' \rangle \geq 1$ for $\bd' \in \bP$, hence obviously $\diff
(\bd', \bd'', [X]) \geq 0$ for all $(\bd', \bd'', [X]) \in \fZ_p$
if $p \geq n$. Moreover, this inequality is strict if $\langle
\bd', \bh \rangle > 1$ and $p > n$.
\end{proof}

\subsection{} 
We consider now the tubular case.

\begin{prop} \label{prop_tub}
Let $\delta = 0$. If $p \geq n + 1$ then $\calZ (p \bh)$ is a set
theoretic complete intersection and $k [\mod_\Lambda (\bd)]$ is a
free $\SI [\mod_\Lambda (\bd)]$-module. The number of irreducible
components of $\calZ (p \bh)$ is
\[
(p - n) m_1 \cdots m_n + \sum_{l \in [1, n - 1]} \sum_{i_1 <
\cdots < i_l \in [1, n]} m_{i_1} \cdots m_{i_l} + 1
\]
for $p > n + 1$.
\end{prop}

\begin{proof}
Fix $(\bd', \bd'', [X]) \in \fZ_p$. Observe that if $\langle \bd',
\bh \rangle = 1$ then $\langle \bd', \bd' \rangle = 1$, according
to Lemma~\ref{lemmdP}, and
\[
\diff (\bd', \bd'', [X]) = p - q \geq 0.
\]
On the other hand, if $\langle \bd', \bh \rangle > 1$ then it
follows from Lemma~\ref{lemmineq} that $\langle \bd', \bd' \rangle
\geq 0$, hence
\[
\diff (\bd', \bd'', [X]) \geq (p - n) (\langle \bd', \bh \rangle -
1) - 1 \geq 0
\]
provided $p \geq n + 1$. Moreover, this inequality is strict if $p
> n + 1$. Now the claim follows again from
Proposition~\ref{propcompbis}, Lemma~\ref{lemmirr} and
Theorem~\ref{theofirst}.
\end{proof}

\subsection{} 
It remains to consider the case $0 < \delta < 1$. We start with
the following observation.

\begin{lemm}
If $0 < \delta < 1$ then $4 \delta + n + 1 < \ffrac{1}{1 - \delta}
(n + 1)$.
\end{lemm}

\begin{proof}
Recall that $n \geq 3$. Consequently,
\[
\ffrac{1}{1 - \delta} (n + 1) = \ffrac{n + 1}{1 - \delta} \delta +
(n + 1) > 4 \delta + n + 1
\]
and the claim follows.
\end{proof}

\subsection{} 
For a fixed $p \geq 1$ consider the real-valued function $f$ given
by
\[
f (t) = - \delta t^2 + t (p - n) + (n - p - 1).
\]

\begin{lemm}
Let $p \geq 1$ and $f$ be as above. If $p \geq \ffrac{1}{1 -
\delta} (n + 1)$ then $f (t) > 0$ for all $t \in [2, p]$.
\end{lemm}

\begin{proof}
It is enough to show that $f (2) > 0$ and $f (p) > 0$. The first
inequality follows from the previous lemma, since
\[
f (2) = p - (4 \delta + n + 1).
\]
On the other hand,
\[
f (p) = p ((1 - \delta) p - (n + 1)) + (n - 1) \geq n - 1 > 0,
\]
which finishes the proof.
\end{proof}

\subsection{} 
Now we can finish our proof.

\begin{prop} \label{prop_wild}
Let $0 < \delta < 1$. If $p \geq \ffrac{1}{1 - \delta} (n + 1)$
then $\calZ (p \bh)$ is a set theoretic complete intersection with
\[
(p - n) m_1 \cdots m_n + \sum_{l \in [1, n - 1]} \sum_{i_1 <
\cdots < i_l \in [1, n]} m_{i_1} \cdots m_{i_l} + 1
\]
irreducible components and $k [\mod_\Lambda (\bd)]$ is a free $\SI
[\mod_\Lambda (\bd)]$-module.
\end{prop}

\begin{proof}
Fix $(\bd', \bd'', [X]) \in \fZ_p$. Observe again that if $\langle
\bd', \bh \rangle = 1$ then $\langle \bd', \bd' \rangle = 1$,
according to Lemma~\ref{lemmdP}, and
\[
\diff (\bd', \bd'', [X]) = p - q \geq 0.
\]
On the other hand, if $t = \langle \bd', \bh \rangle > 1$ then it
follows from Lemma~\ref{lemmineq} that $\langle \bd', \bd' \rangle
\geq - \delta t^2$, hence we obtain from the previous lemma that
\[
\diff (\bd', \bd'', [X]) \geq f (t) > 0
\]
provided $p \geq \ffrac{1}{1 - \delta} (n + 1)$. The claim follows
again from Proposition~\ref{propcompbis}, Lemma~\ref{lemmirr} and
Theorem~\ref{theofirst}.
\end{proof}

\end{document}